\documentclass[fleqn]{article}

\usepackage{arxiv}
\title{The binomial transform of p-recursive sequences and the dilogarithm function}
\author{Stephanie L. Harshbarger\thanks{Work done while an undergraduate student at the University of Nebraska at Kearney. }  \,   and  Barton L. Willis\thanks{Corresponding author (\texttt{willisb@unk.edu}).} \\
Department of Mathematics and Statistics \\
University of Nebraska at Kearney \\
Kearney, Nebraska 68849}

\usepackage[english]{babel}
\usepackage[]{microtype}
\usepackage{enumerate}
\usepackage{comment}
\usepackage[noadjust]{cite}
\usepackage{amsmath, amssymb, amsthm} 
\usepackage{graphicx}
\usepackage{hyperref} 

\newcommand{\BinOp}[2]{\mathrm{B}^{\left (#1,#2 \right)}} 
\newcommand{\MultOp}{\mathrm{M}}
\newcommand{\euler}{\mathrm{e}}
\newcommand{\bigoh}[1]{\ensuremath{\operatorname{O}\left (#1 \right )}}

\newcommand{\Fpq}{{}_{3} \! \operatorname{F}_{2} \!}

\newcommand{\nth}{$\text{n}^\text{th}\,$} 

\DeclareMathOperator{\F}{F}   

\DeclareMathOperator{\ident}{\mathrm{I}}
\DeclareMathOperator{\Sop}{\mathrm{S}}
\DeclareMathOperator{\Li2}{\mathrm{Li}_2}

\newmuskip\pFqmuskip

\newcommand*\pFq[6][8]{%
  \begingroup 
  \pFqmuskip=#1mu\relax
  \begingroup\lccode`\~=`\,
  \lowercase{\endgroup\let~}\pFqcomma
      {}_{#2}\F_{#3}{\left[\genfrac..{0pt}{}{#4}{#5};#6\right]}
  \endgroup
}
\newcommand{\pFqcomma}{\mskip \pFqmuskip}
\newcommand{\imag}{\mathrm{i}\mkern1mu}
\newcommand{\reals}{\mathbf{R}} 
\newcommand{\complex}{\mathbf{C}} 
\newcommand{\integers}{\mathbf{Z}}  
\newcommand{\realpart}{\mathrm{Re}}

\begin{document}

\maketitle

\begin{abstract}
    Using a generalized binomial transform and a novel binomial coefficient identity, we will show that the set of p-recursive sequences is closed under the binomial transform. Using these results, we will derive a new  series representation for the dilogarithm function that converges on \(\complex \setminus [1,\infty)\). Finally, we will show that this series representation results in a scheme for numerical evaluation of the dilogarithm function that is accurate, efficient, and stable.
\end{abstract}

 
\section{Introduction} 

The binomial transform is useful in several contexts, including analytic continuation and series 
acceleration \cite{doi:10.1080/10652469.2016.1231674,Hirofumi}. With an eye toward these applications, we will first show how to derive a generalization of the binomial transform; second, we will show that the set of p-recursive sequences (sequences that satisfy a linear recursion with polynomial coefficients) is closed under the binomial transform; and finally, we will apply these methods to the  Maclaurin series for the dilogarithm function. The result is a series that gives a scheme for numerical evaluation of the dilogarithm function that is accurate, efficient, and stable.

Our notation is fairly standard. The sets of integers, real numbers, and complex numbers are \(\integers, \reals\), and \(\complex\), respectively. We will use subscript modifiers on each of these sets to indicate  various subsets; for example, \(\integers_{\geq 0}\) is the set of nonnegative integers. The real part of a complex number \(x\) is denoted by \(\realpart(x)\), the imaginary unit is \(\imag\) (not \(i\)),  and we use an overline for the complex conjugate. Finally, the identity operator is \(\ident\) and a superscript \(\star\) denotes the operator adjoint.

\section{Generalized binomial transform} 
 
For analytic continuation and series  acceleration, the utility of the binomial transform stems from that fact that it can be derived from a sequence of extrapolated sequences. To show this, we start by defining the backward shift operator \(\Sop\) as
\begin{equation}
    \Sop F = n  \in \integers_{\geq 0} \mapsto \begin{cases} 0 & n = 0 \\ F_{n-1} & n > 0 \end{cases}. 
\end{equation}
 For a sequence \(F^{(0)}\), we define extrapolated sequences \(F^{(1)}, F^{(2)}, \dotsc\) by \(F^{(k)} = \left (\alpha \ident + \beta \Sop \right)^k F^{(0)}\), where \(\alpha, \beta \in \complex\). We call these extrapolated sequences because assuming  \(F^{(0)}\) converges linearly to \(L\), there is a choice of \(\alpha\) and \(\beta\) that makes the convergence of \(F^{(\ell)} \) faster with larger \(\ell\). Specifically, suppose \(F^{(0)} - L \in \bigoh{k \mapsto g^k / k^\mu}\), where \(g, \mu \in \complex\) and \(|g| < 1\). Choosing   \(\alpha = 1/(1-g)\) and \(\beta = g/(g-1)\), we have  
 \(F^{(\ell)} - L \in \bigoh{k \mapsto g^k / k^{\ell + \mu}}\). Although each of these sequences shares the same linear convergence rate \(g\),  the factor  of \(1/k^{\mu+\ell}\) in the leading term of the asymptotic form for \(F^{(\ell)} \) makes convergence of \(F^{(\ell)} \) faster with larger \(\ell\).

Extracting the \nth term of the \nth extrapolated sequence yields a sequence \(n \mapsto  F^{(n)}_n\). We call this sequence the \emph{generalized binomial transform} of \(F^{(0)} \). Defined this way, a calculation shows that the binomial transform operator \(\BinOp{\alpha}{\beta} \) is
 \begin{equation}
   \BinOp{\alpha}{\beta} F=
   n \in \integers_{\geq 0} \mapsto  \sum_{k=0}^n  \binom{n}{k} \alpha^{n-k} \beta^{k} F_k.
 \end{equation}
 In a different context, Prodinger \cite{Prodinger} introduced this form of the binomial transform. We note that depending on the author \cite{oeis,Knuth:1997:ACP:270146}, the standard binomial transformation is either \(\BinOp{1}{1}\) or \(\BinOp{-1}{1}\).
 
 The composition rule for the generalized binomial transform is
\(
  \BinOp{\alpha}{\beta} \BinOp{\alpha^\prime}{\beta^\prime} =
    \BinOp{\alpha + \alpha^\prime  \beta}{\beta \beta^\prime}.
\)
Since \(\BinOp{0}{1}\) is the identity operator, it follows from the composition rule that  for \(\beta \neq 0\), the operator \(\BinOp{\alpha}{\beta}\) is invertible and its inverse is 
\(
    \BinOp{-\alpha/\beta}{1/\beta} 
\). 
Specializing the composition rule to \(\beta = 1\) and \(\beta^\prime =1\), we have \(\BinOp{\alpha}{1} \BinOp{\alpha^\prime}{1} =
\BinOp{\alpha + \alpha^\prime}{1}\). Thus \( \BinOp{\alpha}{1}\) is the \(\alpha\)-fold composition of \(\BinOp{1}{1}\) with itself.

The adjoint of the binomial transform is
\begin{equation}
    {\BinOp{\alpha}{\beta}}^\star F = n \mapsto \sum_{k = n}^\infty  \binom{k}{n}  \alpha^{k-n}  \beta^n F_k.
\end{equation}
 Our interest in the adjoint is the formal identity
 \begin{equation}
    \sum_{k=0}^\infty F_k G_k = \sum_{k=0}^\infty  \big( {\BinOp{-\alpha / \beta}{1/\beta}}^\star G \big)_k   \big(\BinOp{\alpha}{\beta} F \big)_k.  
 \end{equation} 
 We say this is a formal identity because it is valid provided both series converge.  Specializing  to \(G_k =1\) and assuming \(\alpha+\beta \neq 0\) and \(\beta \neq 0\) gives
 \begin{equation}
     \sum_{k=0}^\infty F_k = \sum_{k=0}^\infty    \frac{\beta}{(\alpha+\beta)^{k+1}} \big(\BinOp{\alpha}{\beta} F \big)_k.  
 \end{equation}
 Simplifying the summand shows that  it is a function of the quotient \(\beta/\alpha\) and it does not  depend  individually on \( \alpha \) and \(\beta\). Thus we can assume that \(\beta = 1\). Our identity is an extension of the Euler transform. For a description of the Euler transform, see \S3.9 of the \emph{NIST Handbook of Mathematical Functions}  \cite{NIST:DLMF}.
 
 We will use this   summation  identity to derive a new series representation for the dilogarithm function \(\Li2\). The key to deriving this result is a new binomial coefficient identity.
 
 \section{Binomial coefficient identities}
 
 A sequence that satisfies a linear homogeneous recursion relation with polynomial coefficients is said to be \emph{p-recursive}  \cite{Schneider:2013:CAQ:2541763}. The set of p-recursive sequences is known to be closed under addition and multiplication \cite{Kauers:2011:CT:1993886.1993892, ZEILBERGER1990321}.  We will show that the set of p-recursive sequences is closed under the binomial transform. The key to this result is the binomial coefficient identity
\begin{equation}
     k \binom{n}{k} = n \binom{n}{k} - n \binom{n-1}{k}.
\end{equation}
The proof is a calculation that uses only  simplification and  the factorial representation for the binomial coefficients. Extending this identity by multiplying it by \(k\) and iterating,  allows us to express  \(k^p  \binom{n}{k} \), where \(k \in \integers_{\geq 0}\),  as a linear combination of \(\left \{\binom{n}{k},\binom{n-1}{k}, \binom{n-2}{k}, \dotsc, \binom{n-p}{k} \right \}\) with coefficients that involve only \(n\). Table \ref{BI} displays these results for \(p \) up to three.

Introducing a multiplication operator \(\MultOp\) on the set of sequences defined by \(\MultOp F = n \mapsto n F_n\) and using the identity  \( k \binom{n}{k} = n \binom{n}{k} - n \binom{n-1}{k}\), we can show that
\(
 \BinOp{\alpha}{\beta} \MultOp = \MultOp \big( \ident - \alpha \Sop\big)  \BinOp{\alpha}{\beta}
\), where \(\Sop\) is the backward shift operator.
Consequently, for all \(p \in \integers_{\geq 0}\), we have
\begin{equation}
 \BinOp{\alpha}{\beta} \MultOp^p = \big ( \MultOp (\ident - \alpha S)\big)^p  \BinOp{\alpha}{\beta}.
 \end{equation}
Further, using the  Pascal identity
\(
    \binom{n+1}{k} = \binom{n}{k} + \binom{n}{k-1}
\), we can show that  
\(
    \beta \BinOp{\alpha}{\beta} \Sop^\star = \left( \Sop^\star  - \alpha \ident \right)   \BinOp{\alpha}{\beta}
\). Extending this result to any positive integer power \(p\) of \(\Sop^\star\) yields
\begin{equation}
    \beta^p \BinOp{\alpha}{\beta} \Sop^{\star \, p} = \left(\Sop^{\star} -  \alpha \ident \right)^p \BinOp{\alpha}{\beta}.
\end{equation}
Using these two results, we can express the binomial transform of \(n \mapsto n^p F_{n+q} \) in terms of \(\BinOp{\alpha}{\beta} F\) for all positive integers \(p\) and \(q\). Consequently, we have shown that the set of p-recursive sequences is closed under the generalized binomial transform. 

{\renewcommand{\arraystretch}{1.5}%
\begin{table}[ht]
\centering
\begin{tabular}[p]{| l | l | l | l | l |} \hline
   &
  \(\binom{n-3}{k} \) 
  & \(\binom{n-2}{k} \) 
  & \(\binom{n-1}{k} \) 
  &  \(\binom{n}{k} \)  \\ \hline  \hline
  
\(\binom{n}{k} \)    
& \(0\) 
& \(0\) 
& \(0\) 
& \(1\) \\
\(k \, \binom{n}{k} \)  
& \(0\)
& \(0\) 
& \(-n\) 
& \(n\) \\  
\(k^2 \binom{n}{k} \)  
& \(0\) 
& \(\left( n-1\right) n \) 
& \(-n\,\left( 2n-1\right) \) 
& \( n^2 \) \\ 
\(k^3 \binom{n}{k} \)  
& \(-\left( n-2\right) \,\left( n-1\right) n\) 
&  \(3{{\left( n-1\right) }^{2}}n \) 
& \(-n\,\left( 3{{n}^{2}}-3n+1\right) \)  & \({{n}^{3}} \) \\ \hline
\end{tabular}
\caption{Each row of this table expresses 
\(k^p \binom{n}{k} \) as a linear combination of 
\(\{\binom{n}{k},\binom{n-1}{k}, \binom{n-2}{k}, \dotsc \binom{n-p}{k}\}\) where the coefficients are functions of \(n\) only.  The third row, for example, corresponds to the identity 
\(k^2 \binom{n}{k} =n^2  \binom{n}{k}- n \left(2
 n-1\right) \binom{n-1}{k}  +  \left(n-1\right)\,n \binom{n-2}{k} . \) 
}\label{BI}
\end{table}
}

\section{The dilogarithm function}

The dilogarithm function \(\Li2\) can be defined by its Maclaurin series \cite{NIST:DLMF}
\begin{equation}
    \Li2(x) = \sum_{k=0}^\infty \frac{x^{k+1}}{(k+1)^2}.
\end{equation}
Inside the unit circle, the series converges linearly; on the unit circle, it converges sublinearly, and outside the unit circle, it diverges. 
The summand of the Maclaurin series, call it \(Q\), is p-recursive.  Thus we consider the convergence set for the formal identity
\(
    \Li2(x) = \sum_{k=0}^\infty \widehat Q_k / (\alpha+1)^{k+1},
\)
where \(\widehat{Q} = \BinOp{\alpha}{1} Q\). 

The sequence \(Q\) satisfies the recursion
\(
    \left( k+2\right)^2 Q_{k+1} =  \left( k+1\right)^2 Q_{k}
\).
Using Table \ref{BI}, the recursion for \(\widehat{Q}\) has the form \(0=P_0(n)  \widehat{Q}_n + P_1(n)  \widehat{Q}_{n+1}  
P_2(n)  \widehat{Q}_{n+2} + P_3(n)  \widehat{Q}_{n+3}\), where the polynomials \(P_0\) through \(P_3\) are
\begin{align}
    P_0(n) &= -\alpha^2 (\alpha+x) (n+1)(n+2), \\
    P_1(n) &= \alpha(n+2)(3n\alpha+2n x+8\alpha + 5 x), \\
    P_2(n) &= - \left(3\,\alpha+x\right) n^2\ -\left(19\,\alpha+6\,x\right) n\,-26\,\alpha-9 \, x, \\
    P_3(n) &= (n+2)(n+6).
\end{align}
Assuming \( \alpha \neq -x\), a fundamental solution set for this recursion is 
\begin{equation}
    \left \{  n \mapsto \frac{\alpha^n}{n+1}, \quad   n \mapsto \frac{\alpha^n}{n+1} \sum_{k=0}^n \frac{1}{k+1}, \quad
    n \mapsto  \frac{(\alpha+x)^n}{n^2} \left(1 +  \bigoh{1/n}   \right)
    \right \}.
\end{equation}
The first two members of this set are exact, but the third is an asymptotic solution that is valid toward infinity.

Both  \(\alpha = -x\) and \(\alpha = 0\) are special cases. For \(\alpha = -x\), the order of the recursion is reduced from three to two. For this case, one solution to the recursion is \(n \mapsto (-x)^n /(1+n) \). Since this series diverges everywhere outside the unit circle, we will discard it. Similarly, the case \(\alpha = 0\) is not pertinent.

 The fundamental solution set shows that  the formal series converges linearly,  provided that
\begin{equation}
     \max \left( \left| \frac{\alpha}{\alpha+1} \right|, \left| \frac{\alpha+ x}{\alpha + 1} \right| \right) < 1,
     \mbox{ and } \alpha \in \complex_{\neq -x, \neq -1}.
\end{equation}
The convergence set is maximized when \(\left| \frac{\alpha}{\alpha+1} \right| =  \left| \frac{\alpha+ x}{\alpha + 1} \right|\). Assuming \(x \in \reals\), the convergence set is maximized when \(\alpha=-x/2\). For this choice, the  linear convergence rate is \(|x/(x-2)| \) and the series converges in the half plane  \(\realpart(x) < 1\). 
For \(x \in \complex \setminus [1, \infty) \), the convergence set is maximized when
\begin{equation}
    \alpha =  \frac{\euler^{\imag \theta}}{\euler^{\imag \theta} - 1} x \mbox{ , where  } 
    \euler^{\imag \theta}  = \pm  \sqrt{\frac{\overline{x}-1}{x-1}}.
\end{equation}
Setting \(x = 1 + R \exp( \imag \omega) \), where \(R \in \reals_{\geq 0} \) and \(\omega \in [0, 2 \pi) \), the minimum of the linear convergence rate \(|\alpha/(\alpha+1)| \) is
\begin{equation}
    \min\left( \frac{2 R \cos{(\omega)}+{{R}^{2}}+1}{{{\left( R-1\right) }^{2}}},\frac{2 R \cos{(\omega)}+{{R}^{2}}+1}{{{\left( R+1\right) }^{2}}}\right).
\end{equation}
For \(\omega \in (0, 2 \pi)\), or equivalently for \(x \in \complex \setminus [1,\infty)\), the linear convergence rate is less than one. Consequently, there is a value of \(\alpha\) that makes the series \(\sum_{k=0}^\infty \widehat Q_k / (\alpha+1)^{k+1} \) converge on \(\complex \setminus [1,\infty)\). Although this is a satisfying result, its additional complexity over the choice \(\alpha = -x/2\) is  erased by the fact that \(\Li2\) satisfies several functional identities that allow the convergence set \(\realpart(x) < 1 \) to be adequate, at least for numerical evaluation.

Returning to the choice \(\alpha = -x/2\), the recursion for the entire  summand \(W_k = \widehat{Q}_k / (1+\alpha)^{k+1} \) is
\begin{equation}
   (x-2)^3 (n+4)^2 W_{n+3} =  -x^3 (n+1)(n+2)W_n 
   +x^2 (x-2) (n+2)^2 W_{n+1} + x (x-2)^2 (n+3)(n+4) W_{n+2}.
\end{equation}
In terms of the forward shift operator \(S^\star\), the recursion relation factors as
\begin{equation}
    \big ( (x-2)^2 (n+4) S^{\star 2} - x^2 (n+2) \big)
    \big ((x-2) (n+2) S^\star - x (n + 1) \big) W_k = 0.
\end{equation}
The three initial values of the sequence \(W\) are
\begin{align}
    W_0 &=  x/(1-x/2), \\ 
    W_1 &= -x^2/\left(4(1-x/2)^2 \right), \\ 
	W_2 &= x^3/\left(9(1-x/2)^3 \right).
\end{align}
For a series that converges in the half-plane \(\realpart(x) < 1/2\), see \cite{MaxieSchmidt}.

In the next section, we will investigate the practical considerations of using this series to numerically evaluate \(\Li2\).

\section{Accuracy, efficiency, and stability}

For our series representation to be useful for numerical evaluation, the sum must be well conditioned (accuracy), the
convergence must be fast (efficiency), and every solution to the fundamental solution set to the recursion for the summand must converge to zero (stability).

Of these three conditions, we have already shown  that each member of the fundamental solution set to the recursion relation converges to zero when \(\realpart(x) < 1\). Thus the recursion for \(W\) is stable. 

\begin{figure}[ht!]
\includegraphics[width=0.5\textwidth]{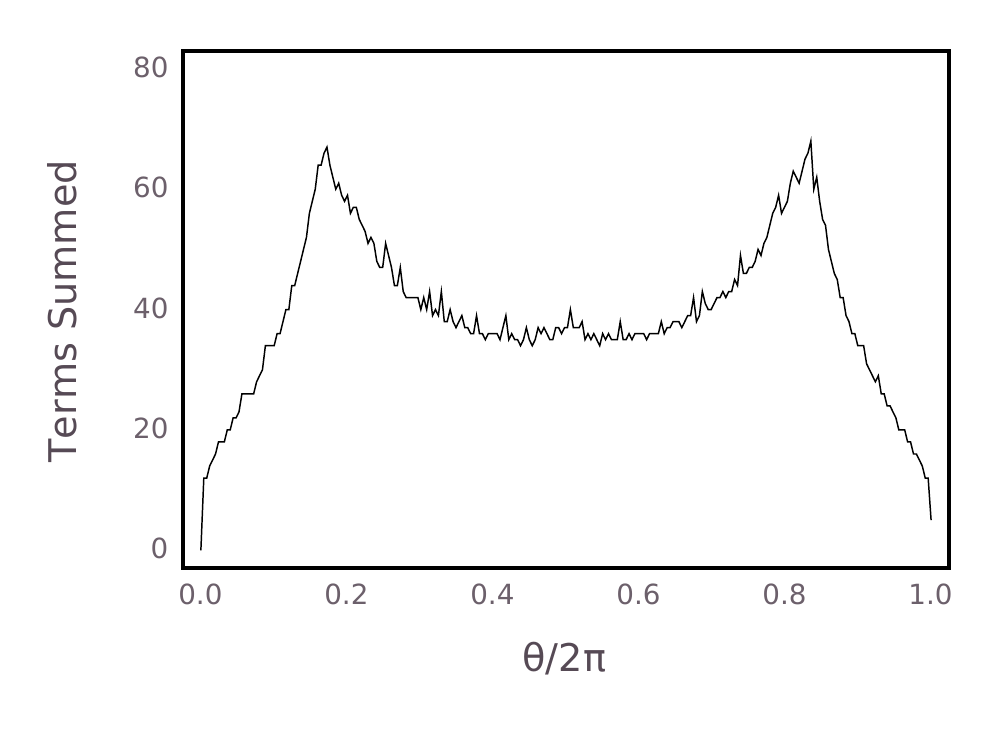}
\centering
\caption{This graph shows the number of terms that need to be summed to achieve full accuracy using IEEE binary64 numbers to compute the value of \(\Li2(\exp( \imag \theta /  2 \pi))  \) for \(\theta \in [0, 2 \pi]\). The maximum number of terms is less than 70.}\label{nbr_summed}
\end{figure}

We can achieve greater efficiency by leveraging two functional identities. The algorithm can  automatically choose between them to minimize the linear convergence rate of \(|x/(2-x) |\). These functional identities are (see \cite{NIST:DLMF})
\begin{align}
    \Li2(x) + \Li2(1-x)  &= \pi^{2}/6 - \ln(x) \, \ln (1-x),  \quad x \in \complex_{\neq 0, \neq 1}, \\
    \mathrm{Li}_{2}\left(x\right)+\mathrm{Li}_{2}\left(\frac{1}{x}\right) & =-\pi^2 /6 - (\ln\left(-x\right))^{2}/2, \quad x \in \complex \setminus [0,1].
\end{align}
The second identity is the  \(\Li2\) reciprocal formula. Choosing between these identities and using our series to compute \(\Li2\) on the unit circle, the largest number of terms that must be summed to achieve full accuracy using IEEE binary64 numbers is less than 70; see Figure \ref{nbr_summed}.  Inside the unit circle, the number of terms needed is less than 70; and outside, the \(\Li2\) reciprocal formula reduces this evaluation to a number on the inside of the unit circle.

\begin{figure}[ht!]
\includegraphics[width=0.5\textwidth]{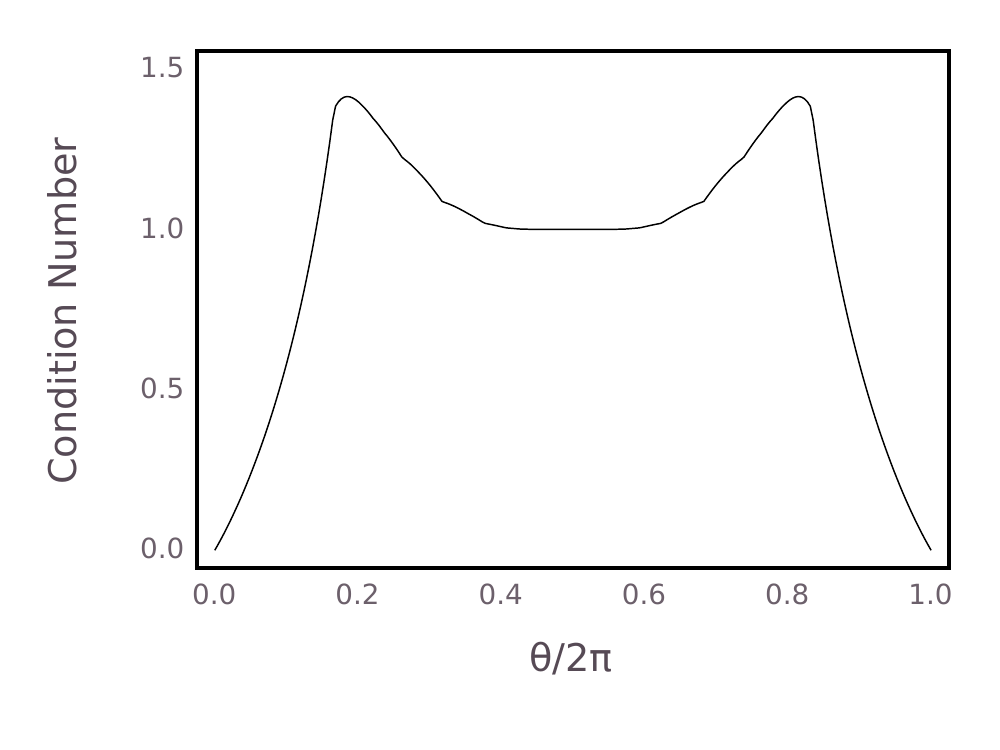}
\centering
\caption{The condition number for summing \(\sum_k W_k(\exp(\imag \theta)) \). On the unit circle, the condition number is apparently bounded
above by \(3/2\).}\label{cndX}
\end{figure}
Finally, we study the condition number of the sum.  Recall that the condition number  \cite{Higham:2002:ASN} of a sum \(\sum_k  W_k\) is the quotient \(\sum_k |W_k| / | \sum_k W_k | \). Using Kahan summation \cite{Kahan:1965:PRR:363707.363723}, the floating point rounding error is bounded by the machine epsilon (\(2^{-53}\) for a IEEE binary64 number) times twice the condition number. Again, automatically choosing between the functional identities, the condition number for the sum \(\sum_k W_k\) is shown in Figure \ref{cndX} for inputs on the unit circle. The condition number is apparently bounded above by \(3/2\); thus the sum is well conditioned on the unit circle.

For testing, we implemented the algorithm in the Julia language \cite{bezanson2017julia}. Our implementation uses Kahan summation and it accumulates the condition number for the sum. The condition number indicates the total rounding error; for details, see Higham \cite{Higham:2002:ASN}. Finally, the method is  generic for both real and complex IEEE floats, as well as real and complex extended precision floating point numbers.

\section{Acknowledgments}
The work by Stephanie Harshbarger was supported by the  University of Nebraska at Kearney  Undergraduate Research Fellows Program. We used the Maxima computer algebra system  \cite{maxima} to do the calculations in this paper. We thank the volunteers who make Maxima freely available.
\bibliography{precursive}
\bibliographystyle{plain}
\end{document}